\documentclass{amsart}
\usepackage{amssymb}
\def\version{Version W06v6c.tex last changed 05 November 2002}
\usepackage{a4}
\begin{document}
\date{\version}
\newtheorem{theorem}{Theorem}[section]
\newtheorem{lemma}[theorem]{Lemma}
\newtheorem{remark}[theorem]{Remark}
\newtheorem{definition}[theorem]{Definition}
\newtheorem{corollary}[theorem]{Corollary}
\newtheorem{example}[theorem]{Example}
\font\pbglie=eufm10
\def\SS{\mathcal{S}}
\def\BB{\mathcal{B}}
\def\RR{\mathcal{R}}
\def\Gr{\operatorname{Gr}}
\def\Range{\operatorname{Range}}
\def\Ker{\operatorname{Ker}}
\def\Span{\operatorname{Span}}
\def\Tr{\operatorname{Trace}}
\makeatletter
  \renewcommand{\theequation}{%
   \thesection.\alph{equation}}
  \@addtoreset{equation}{section}
 \makeatother
\date{\version}
\title[Szabo manfiolds]
{Nilpotent Szab\'o, Osserman and Ivanov-Petrova pseudo-Riemannian manifolds}
\author{B. Fiedler and P. Gilkey}
\begin{address}{BF: Mathematics Institute, University of Leipzig, Augustusplatz 10/11,
04109 Leipzig,
Germany}\end{address}
\begin{email}{bernd.fiedler.roschstr.leipzig@t-online.de}\end{email}
\begin{address}{PG: Max Planck Institute for Mathematics in the Sciences, 
Inselstrasse 22-26, 04103 Leipzig, Germany and Mathematics Dept., University of Oregon,
Eugene Or 97403 USA}
\end{address}
\begin{email}{gilkey@darkwing.uoregon.edu}\end{email}
\begin{abstract} We exhibit pseudo Riemannian manifolds which are Szab\'o nilpotent of
arbitrary
order, or which are Osserman nilpotent of arbitrary
order, or which are Ivanov-Petrova nilpotent of order 3.\end{abstract}
\maketitle
\section{Introduction}
Let $R$ be the Riemann curvature tensor of a pseudo-Riemannian manifold $(M,g)$ of
signature $(p,q)$. The {\it Szab\'o operator} $\SS$ is the self-adjoint linear map which
is
characterized by the identity:
$$g(\SS(x)y,z)=\nabla R(y,x,x,z;x).$$
One says that $(M,g)$ is {\it Szab\'o} if the eigenvalues of $\SS(x)$ are constant on the
pseudo-spheres of
unit timelike and spacelike vectors:
$$S^\pm(M,g):=\{x\in TM:g(x,x)=\pm1\}.$$

Szab\'o \cite{Sz91} used techniques from algebraic topology to show in the Riemannian
setting ($p=0$) that any such metric
is locally symmetric. He used this observation to give a simple proof that any $2$ point
homogeneous space is either flat or is a rank $1$ symmetric space. Subsequently Gilkey
and Stavrov \cite{GS02} extended
his results to show that any Szab\'o Lorentzian $(p=1)$ manifold has constant sectional
curvature. By replacing $g$ by $-g$,
one can interchange the roles of $p$ and of $q$, thus these results apply to the cases
$q=0$ and $q=1$ as well.

The eigenvalue zero is distinguished. One says that $(M,g)$ is {\it Szab\'o nilpotent of
order $n$} if $\SS(x)^n=0$ for
every $x\in TM$ and if there exists a point $P_0\in M$ and a tangent vector $x_0\in
T_{P_0}M$ so that $\SS(x_0)^{n-1}\ne0$.
One says that $(M,g)$ is {\it Szab\'o nilpotent} if $(M,g)$ is Szab\'o nilpotent of order
$n$ for some $n$. Note that
$(M,g)$ is Szab\'o nilpotent if and only if $0$ is the only eigenvalue of $\SS$;
consequently any Szab\'o nilpotent manifold
is Szab\'o. There is some evidence 
\cite{GIS02,St03} to suggest, conversely, that any Szab\'o manifold is Szab\'o nilpotent.

If $(M,g)$ is Szab\'o nilpotent  of order $1$, then
$\SS(x)=0$ for all $x\in TM$. This implies \cite{GS02} that $\nabla R=0$ so $(M,g)$ is a
local symmetric space; this is to be
regarded, therefore, as a trivial case.
 Gilkey, Ivanova, and Zhang
\cite{GIZ02} have constructed pseudo-Riemannian manifolds of any signature $(p,q)$ with
$p\ge2$ and $q\ge2$ which are
Szab\'o nilpotent of order $2$; these were the only previously known examples of Szab\'o
manifolds which were not local
symmetric spaces. In this brief note, we shall construct pseudo-Riemannian metrics
$g_{n}$ on $\mathbb{R}^{n+2}$ which are
Szab\'o nilpotent of order $n\ge2$; the metric will be {\it balanced} (i.e. $p=q$) if $n$
is even and {\it almost balanced}
(i.e. $p=q\pm1$) if $n$ is odd. By taking an isometric product with a suitable flat
manifold, the signature can be increased
without changing the order of nilpotency.

\begin{theorem}\label{thm1.1} Let $n\ge2$. There exists a pseudo-Riemannian metric
$g_{n}$ on $\mathbb{R}^{n+2}$
which is Szab\'o nilpotent of order $n$. If $n=2p$, then $g_{n}$ has signature
$(p+1,p+1)$; if $n=2p+1$, then $g_{n}$
has signature $(p+1,p+2)$.
\end{theorem}

The {\it Jacobi operator} is defined analogously; it is characterized by the identity:
$$g(J(x)y,z)=R(y,x,x,z).$$
One says that $(M,g)$ is {\it Osserman} if the eigenvalues of $J$ are constant on
$S^\pm(M)$. In the Riemannian setting,
Osserman wondered \cite{Oss90} if this implied $(M,g)$ was a 2 point homogeneous space.
This
question has been answered in the affirmative in the Riemannian setting \cite{C88,N02}
for dimensions $\ne8,16$, and in
all dimensions in the Lorentzian setting \cite{BBG97,GKV97}.

We shall say that $(M,g)$ is {\it Osserman nilpotent of order n} if $J(x)^n=0$ for every
$x\in TM$ and if there exists a point $P_0\in M$ and a tangent vector $x_0\in T_{P_0}M$
so that $J(x_0)^{n-1}\ne0$, i.e. $0$
is the only eigenvalue of $J$. Such manifolds are necessarily Osserman. Osserman
nilpotent manifolds of orders $2$ and $3$
have been constructed previously
\cite{BBGR97,GVV02,GVV98,G01}. These manifolds need not be homogeneous, thus the question
Osserman raised has a negative
answer in the higher signature setting. A byproduct of our investigation of Szab\'o
manifolds yields new examples of Osserman
manifolds; again, the signature can be increased by taking isometric products with flat
factors.

\begin{theorem}\label{thm1.2} Let $n\ge2$. There exists a pseudo-Riemannian metric
$\tilde g_{n}$ on $\mathbb{R}^{n+2}$
which is Osserman nilpotent of order $n$. If $n=2p$, $\tilde g_{n}$ has signature
$(p+1,p+1)$; if $n=2p+1$,
$\tilde g_{n}$ has signature $(p+1,p+2)$.
\end{theorem}

If $\{f_1,f_2\}$ is an oriented orthonormal basis for a non-degenerate oriented $2$ plane
$\pi$, we define the skew-symmetric
curvature operator by setting $\mathcal{R}(\pi):=R(f_1,f_2)$. We say $(M,g)$ is
Ivanov-Petrova nilpotent of order $n$ if
$\mathcal{R}(\pi)^n=0$ for any non-degenerate oriented $2$ plane $\pi$ and if there
exists $\pi$ so
$\mathcal{R}(\pi)^{n-1}\ne0$. We refer to
\cite{G01} for further details concerning Ivanov-Petrova manifolds. Another byproduct of
our investigation yields new examples
of these manifolds:
\begin{theorem}\label{thm1.3}\ There exist Ivanova-Petrova pseudo-Riemannian manifolds
which are nilpotent of order 2 and of
order 3.\end{theorem}

Here is a brief outline to the paper. In Section \ref{Sect2}, we give a general procedure
for constructing pseudo-Riemannian
manifolds with certain kinds of curvature and covariant derivative curvature tensors. We
apply this procedure in Section
\ref{Sect3} to complete the proof of Theorem
\ref{thm1.1}. Lemma \ref{lem3.1} deals with the cases $n=2$ and $n=3$, Lemma
\ref{lem3.2} deals with the case $n=2\ell+1\ge5$, and Lemma \ref{lem3.3} deals with the
case $n=2\ell+2\ge4$. In
Section \ref{Sect4}, we prove Theorem \ref{thm1.2} and in Section \ref{Sect5}, we prove
Theorem \ref{thm1.3}.

One can also work with the Jordan normal form; one says $(M,g)$ is Jordan Szabo (resp.
Jordan Osserman or Jordan IP) if the
Jordan normal form of $\SS$ (resp. $J$ or $\RR$) is constant on the appropriate domains
of definition. The examples
constructed in this paper do {\bf not} fall into this framework; in particular, there are
no known Jordan Szabo
pseudo-Riemannian manifolds which are not locally symmetric.

\section{A family of pseudo-Riemannian manifolds}\label{Sect2}

We introduce the following notational conventions.
Let $(x,u_1,...,u_\nu,y)$ be coordinates on $\mathbb{R}^{\nu+2}$. We shall use several
different notations for the coordinate
frame:

$$\BB=\{e_0,e_1,...,e_{\nu+1}\}=\{X,U_1,...,U_\nu,Y\}:=\{\partial_x,\partial_{u_1},...,\partial_{u_\nu},\partial_y\}.$$
Let indices $i,j,...$ range from $0$ through $\nu+1$ and index the full coordinate frame.
Let indices $a,b$ range
from $1$ through $\nu$ and index the tangent vectors $\{U_1,....,U_\nu\}$. In the
interests of brevity, we shall give non-zero entries in a metric $g$, curvature tensor
$R$, and covariant derivative curvature
tensor $\nabla R$ up to the obvious
$\mathbb{Z}_2$ symmetries. 

\begin{lemma}\label{lem2.1}
Let
$f=f(u)$ be a smooth function on $\mathbb{R}^\nu$ and let $\Xi$ be a constant invertible
symmetric
$\nu\times\nu$ matrix. Define a metric $g_f$ on
$\mathbb{R}^{\nu+2}$ by setting:
$$\begin{array}{llll}
g_f(X,X)=f(u),&
g_f(X,Y)=1,&\text{and}&
g_f(U_a,U_b)=\Xi_{ab}.\vphantom{\vrule height 12pt}
\end{array}$$
All other scalar products equals zero.
\begin{enumerate}
\item Then the non-zero entries in $R_{g_f}$ are
$R_{g_f}(X,U_a,U_b,X)=-\textstyle\frac12U_aU_b(f)$.
\item The non-zero entries in $\nabla R_{g_f}$ are
$\nabla R_{g_f}(X,U_a,U_b,X;U_c)=
   -\frac12U_aU_bU_c(f)$.
\end{enumerate}
\end{lemma}

\begin{proof}
Since $d\Xi=0$, the non-zero Christoffel symbols of the first kind are:
\begin{equation}\label{eqn2.a}
   \Gamma_{a00}=\Gamma_{0a0}=-\Gamma_{00a}=\textstyle\frac12U_a(f).
\end{equation}
Let $\Xi^{ab}$ be the inverse matrix. We adopt the Einstein convention and sum over
repeated indices to compute:
$$\begin{array}{ll}

\Gamma_{ijb}=g(\nabla_{e_i}e_j,e_b)=g(\Gamma_{ij}{}^ke_k,e_b)=\Gamma_{ij}{}^a\Xi_{ab}&\text{so
}
\Gamma_{ij}{}^a=\Xi^{ab}\Gamma_{ijb},\\

\Gamma_{ij\nu+1}=g(\nabla_{e_i}e_j,e_{\nu+1})=g(\Gamma_{ij}{}^ke_k,e_{\nu+1})=\Gamma_{ij}{}^0&\text{so
}
\Gamma_{ij}{}^0=0,\vphantom{\vrule height 12pt}\\

\Gamma_{ij0}=g(\nabla_{e_i}e_j,e_0)=g(\Gamma_{ij}{}^ke_k,e_0)=f\Gamma_{ij}{}^0+\Gamma_{ij}{}^{\nu+1}&\text{so
}
\Gamma_{ij}{}^{\nu+1}=\Gamma_{ij0}.\vphantom{\vrule height 12pt}
\end{array}$$
Thus the non-zero Christoffel symbols of the second kind are:
\begin{equation}\label{eqn2.b}
\Gamma_{a0}{}^{\nu+1}=\Gamma_{0a}{}^{\nu+1}=\textstyle\frac12U_a(f)\quad\text{and}
\quad\Gamma_{00}{}^a=-\textstyle\frac12\sum_b\Xi^{ab}U_b(f).
\end{equation}
The components of the curvature tensor relative to the coordinate frame are:
\begin{equation}\label{eqn2.c}
R_{ijkl}=e_i\Gamma_{jkl}-e_j\Gamma_{ikl}+\textstyle\sum_n\{\Gamma_{inl}\Gamma_{jk}{}^n
 -\Gamma_{jnl}\Gamma_{ik}{}^n\}.\end{equation}
By equation (\ref{eqn2.b}), $\Gamma_{ik}{}^0=\Gamma_{jk}{}^0=0$. By equation
(\ref{eqn2.a}), 
$\Gamma_{i,\nu+1,k}=\Gamma_{j,\nu+1,k}=0$. Thus the index $n$ in equation (\ref{eqn2.c})
is neither $0$ nor $\nu+1$. Thus by
equation (\ref{eqn2.a}) and equation (\ref{eqn2.b}), $i=j=k=l=0$. This shows that the
terms which are quadratic in $\Gamma$
play no role in equation (\ref{eqn2.c}). Assertion (1) then follows from equation
(\ref{eqn2.a}).

The covariant derivative of the curvature tensor is given by:
\begin{equation}\label{eqn2.d}
 R_{ijkl;n}=e_nR_{ijkl}-\textstyle\sum_p\{\Gamma_{ni}{}^pR_{pjkl}+\Gamma_{nj}{}^pR_{ipkl}
   +\Gamma_{nk}{}^pR_{ijpl}+\Gamma_{nl}{}^pR_{ijkp}\}.\end{equation}
By equation (\ref{eqn2.b}) $\Gamma_{**}{}^0=0$. Thus we may assume $p\ne0$ in equation
(\ref{eqn2.d}). Furthermore, by assertion (1),
$R_{\nu+1***}=R_{*\nu+1**}=R_{**\nu+1*}=R_{***\nu+1}=0$
so we may also assume
$p\ne\nu+1$ in equation  (\ref{eqn2.d}). Thus
$\Gamma_{ni}{}^pR_{pjkl}=0$ unless
$i=j=0$ and similarly $\Gamma_{nj}{}^pR_{ipkl}=0$ unless $i=j=0$. Thus these two terms
cancel. Similarly
$\Gamma_{nk}{}^pR_{ijpl}$ cancels $\Gamma_{nl}{}^pR_{ijkp}$. Thus
$R_{ijkl;n}=e_nR_{ijkl}$ and assertion (2) follows.
\end{proof}

\begin{remark}\rm Let $\rho$ be the associated Ricci tensor;
$\rho(\xi,\xi)=\Tr(J(\xi))$. We have
$\rho(e_i,e_j)=\sum_{kl}g^{kl}R(e_i,e_k,e_l,e_j)$. Since $R$ vanishes
on
$e_{\nu+1}$,  we may sum over $k,l\le\nu$. Since $g^{0k}=g^{k0}=0$ for
$k\le\nu$,
$\rho(e_i,e_j)=\sum_{ab}g^{ab}R(e_i,e_a,e_b,e_j)$. Thus
$\rho(e_i,e_j)=0$ for $(i,j)\ne(0,0)$ and the only non-zero
entry of the Ricci tensor is 
$\rho(e_0,e_0)=-\textstyle\frac12\sum_{ab}\Xi^{ab}\partial_a\partial_bf$. 
The associated Jacobi operator will be nilpotent if and only if 
this sum vanishes. Raising indices yields a Ricci operator $\hat\rho$
with the property that 
$\hat\rho(e_0)=-\textstyle\frac12\sum_{ab}\Xi^{ab}\partial_a\partial_bf$ and 
$\hat\rho(e_i)=0$ for $i>0$. Thus $\hat\rho^2=0$ so the Ricci operator 
is nilpotent of order $2$ and non-trivial if
and only if $g$ is not Osserman.
\end{remark}

If $f$ is quadratic, then $R$  is constant on the coordinate frame; if $f$ is
cubic, then $\nabla R$ is constant on the coordinate frame. However, these tensors are
not curvature homogeneous in the sense
of Kowalski, Tricerri, and Vanhecke
\cite{KTV92} since the metric relative to the coordinate frames is not constant.

The tensors of Lemma \ref{lem2.1} are related to hypersurface theory. Let $M$ be a
non-degenerate
hypersurface in $\mathbb{R}^{(a,b)}$; we assume $M$ is spacelike but similar remarks hold
in the timelike setting. Let $L$ be
the associated second fundamental form and let $S=\nabla L$ be the covariant derivative
of $L$; $L$ is a totally symmetric
$2$ form and $S$ is a totally symmetric
$3$ form. We may then, see for example \cite{G01}, express:
\begin{eqnarray}
R_L(x_1,x_2,x_3,x_4)&=&L(x_1,x_4)L(x_2,x_3)-L(x_1,x_3)L(x_2,x_4),\nonumber\\
\nabla
R_{L,S}(x_1,x_2,x_3,x_4;x_5)&=&S(x_1,x_4,x_5)L(x_2,x_3)+L(x_1,x_4)S(x_2,x_3,x_5)\nonumber\\
&-&S(x_1,x_3,x_5)L(x_2,x_4)-L(x_1,x_3)S(x_2,x_4,x_5).\label{eqn2.e}
\end{eqnarray}

If $L$ is an arbitrary symmetric $2$ tensor and if
$S$ is an arbitrary totally symmetric
$3$ tensor, then we may use equation (\ref{eqn2.e}) to define tensors we continue to
denote by $R_L$ and $\nabla
R_{L,S}$. We refer to
\cite{F01} for the proof of assertion (1) and to
\cite{F02} for the proof of assertion (2) in the following result:
\begin{theorem}\label{thm2.1}\ \begin{enumerate}
\item The tensors $R_L$ which are defined by a symmetric $2$ form $L$ generate the space
of all algebraic curvature tensors.
\item The tensors $\nabla R_{L,S}$ which are defined by a symmetric $2$ form $L$ and by a
totally symmetric $3$ form $S$
generate the space of all algebraic covariant derivative curvature tensors.
\end{enumerate}
\end{theorem}

The tensors of Lemma \ref{lem2.1} (2) are of this form. Let
$$f(u):=-\textstyle\frac13\sum_{a,b,c}c_{a,b,c}u_au_bu_c$$
be a cubic polynomial in the $u$ variables which is independent of $x$ and $y$. Then:
\begin{eqnarray*}
&&\nabla R=\nabla R_{L,S}\ \ \text{for}\ \ 
L(\partial_i,\partial_j):=\delta_{0,i}\delta_{0,j}\ \ \text{and}\ \ 
S(\partial_i,\partial_j,\partial_k):=-\textstyle\frac12\partial_i\partial_j\partial_kf.\end{eqnarray*}

\section{Nilpotent Szab\'o manifolds}\label{Sect3}
In this section we will use Lemma \ref{lem2.1} to prove Theorem \ref{thm1.1} by choosing
$\Xi$ and $f$ appropriately. We shall consider metrics of the form:
$$g(X,X)=f(t,u,v),\ g(X,Y)=1,\ g(T,T)=1,\ g(U_a,V_b)=\delta_{ab};$$
 the spacelike vector $T$ will not be present in some cases. The vectors $\{U_a,V_a\}$
are a hyperbolic pair.

We begin
by discussing the cases $n=2$ and $n=3$.
\begin{lemma}\label{lem3.1}\ \begin{enumerate}
\item Let $\BB_2:=\{X,U,V,Y\}=\{\partial_x,\partial_{u},\partial_{v},\partial_y\}$ be the
coordinate frame
on $\mathbb{R}^4$ relative to the coordinate system $(x,u,v,y)$. Define a 
metric
$g_2$ by:
$$\textstyle g_2(X,X)=-\frac13u^3,\ 
g_2(X,Y)=1,\ g_2(U,V)=1.$$
Then $g_2$ has signature $(2,2)$ on $\mathbb{R}^4$ and $g_2$ is Szab\'o nilpotent of
order $2$.
\smallbreak\item Let
$\BB_3:=\{X,T,U,V,Y\}=\{\partial_x,\partial_t,\partial_u,\partial_v,\partial_y\}$ be the
coordinate
frame on
$\mathbb{R}^5$ relative to the coordinate system $(x,t,u,v,y)$. Define a metric $g_3$ by:
$$g_3(X,X)=\textstyle-tu^2,\ g_3(T,T)=1,\ g_3(U,V)=1, \ g_3(X,Y)=1.$$
Then $g_3$ has signature $(2,3)$ on $\mathbb{R}^5$ and $g_3$ is Szab\'o nilpotent of
order $3$.
\end{enumerate}\end{lemma}

\begin{proof} Let $\BB^*=\{e^0,\ldots,e^{\nu+1}\}$ be the corresponding {\it dual basis}
of $\BB$; it is characterized by
the relations
$g(e_i,e^j)=\delta_i^j$.
For example, we have
\begin{equation}\label{eqn3.x}\BB_2^*=\{Y,V,U,X-fY\}\quad\text{and}\quad
\BB_3^*=\{Y,T,V,U,X-fY\}\end{equation}

By Lemma \ref{lem2.1}, the only non-zero component of $\nabla R_{g_2}$ is
$$\nabla R_{g_2}(X,U,U,X;U)=1.$$
Let $\xi=\xi_0X+\xi_1U+\xi_2V+\xi_3Y$ be a tangent vector. We use equation (\ref{eqn3.x})
to raise indices and conclude:
$$\SS_{g_2}(\xi)X=\xi_1^3 Y-\xi_0\xi_1^2 V,\quad
\SS_{g_2}(\xi)U=-\xi_0\xi_1^2 Y+\xi_0^2\xi_1 V,\quad
\SS_{g_2}(\xi)Y=\SS(\xi)V=0.$$
Thus $\SS_{g_2}(\xi)^2=0$ for all $\xi$ while
$\SS_{g_2}(\xi)\ne0$ for generic $\xi$. Assertion (1) now follows.

Similarly, the only non-zero components of $\nabla R_{g_3}$ are 
$$\nabla R_{g_3}(X,U,U,X;T)=
  \nabla R_{g_3}(X,U,T,X;U)=1.$$
We use equation (\ref{eqn3.x}) to raise indices and compute:
$$\begin{array}{ll}
\SS_{g_3}(\xi)X=\star T+\star Y+\star V,&
\SS_{g_3}(\xi)Y=0,\\
\SS_{g_3}(\xi)T=\star Y+\star V,\\
\SS_{g_3}(\xi)U=\star T+\star Y+\star V,&
\SS_{g_3}(\xi)V=0
\end{array}$$
where $\star=\star(\xi)$ denotes suitably chosen cubic polynomials in the
coefficients of $\xi$ that is generically non-zero; as the precise value of this
coefficient is not important, we
shall suppress it in the interests of notational simplicity. It is now clear that
$\SS_{g_3}(\xi)^3=0$ for all
$\xi$ while
$\SS_{g_3}(\xi)^2$ is generically non-zero.
\end{proof}

Next we consider the case $n=2\ell+1\ge5$. Let
$(x,t,u_2,...,u_{\ell+1},v_2,...,v_{\ell+1},y)$ be coordinates on
$\mathbb{R}^{2\ell+3}$ which define the associated coordinate frame:
$$\BB:=\{X,T,U_2,...,U_{\ell+1},V_2,...,V_{\ell+1},Y\}=
\{\partial_x,\partial_t,\partial_{u_2},...,\partial_{u_{\ell+1}},\partial_{v_2},...,\partial_{v_{\ell+1}},\partial_y\}.$$

\begin{lemma}\label{lem3.2} Let $\ell\ge2$.
Define
a metric $g_{2\ell+1}$ on
$\mathbb{R}^{2\ell+3}$ by setting:
\begin{eqnarray*}
&&g_{2\ell+1}(X,X)=-\textstyle tu_2^2-\textstyle\textstyle\sum_{2\le
a\le\ell}(u_a+v_a)u_{a+1}^2,\\
&&g_{2\ell+1}(X,Y)=1,\quad g_{2\ell+1}(T,T)=1,\quad
   g_{2\ell+1}(U_a,V_b)=\delta_{ab}.\vphantom{\vrule height 12pt}
\end{eqnarray*}
Then $g_{2\ell+1}$ is a metric of signature $(\ell+1,\ell+2)$ and Szab\'o nilpotent of
order $2\ell+1$.
\end{lemma}

\begin{proof} Let $2\le a\le\ell$. By Lemma \ref{lem2.1}, the non-zero components of
$\nabla
R$ are:
$$\begin{array}{lllll}
\nabla R(X,U_2,U_2,X;T)&=&
\nabla R(X,T,U_2,X;U_2)&=&1,\\
\nabla R(X,U_{a+1},U_{a+1},X;U_a)&=&
\nabla R(X,U_{a+1},U_a,X;U_{a+1})&=&1,\\
\nabla R(X,U_{a+1},U_{a+1},X;V_a)&=&
\nabla R(X,U_{a+1},V_a,X;U_{a+1})&=&1.
\end{array}$$
The dual basis is
$\BB^*=\{Y,T,V_2,...,V_{\ell+1},U_2,...,U_{\ell+1},X-fY\}$.
Let $\xi$ be an arbitrary tangent vector. We raise indices and
compute:
$$\begin{array}{lll}
\SS(\xi)X&\in&\Span\{Y,T,U_2,...,U_\ell,
        V_2,...,V_{\ell+1}\},\\
\SS(\xi)Y&=&0,\\
\SS(\xi)T&=&\star Y+\star V_2,\\
\SS(\xi)U_2&=&\star T+\star Y+\star V_2+\star V_3,\\
\SS(\xi)U_a&=&\star U_{a-1}+\star Y+\star V_{a-1}+\star V_a+\star V_{a+1}\quad
\text{for}\quad3\le a\le\ell,\\
\SS(\xi)U_{\ell+1}&=&\star U_\ell+\star Y+\star V_\ell+\star V_{\ell+1}\\
\SS(\xi)V_a&=&\star Y+\star V_{a+1}\quad\text{for}\quad2\le a\le\ell,\\
\SS(\xi)V_{\ell+1}&=&0
\end{array}$$
where $\star $ is a coefficient that is non-zero for generic $\xi$. If $\mathcal{E}$ is a
subspace, let
$\alpha=\beta+\mathcal{E}$ mean that
$\alpha-\beta\in
\mathcal{E}$. We compute:
$$\begin{array}{llll}
\SS(\xi)^\mu U_{\ell+1}&=&\star U_{\ell+1-\mu}+\Span\{V_2,...,V_{\ell+1},Y\},&1\le
\mu\le\ell-1\\
\SS(\xi)^{\ell}U_{\ell+1}&=&\star T+\Span\{V_2,...,V_{\ell+1},Y\},\\
\SS(\xi)^{\mu}U_{\ell+1}&=&\star
V_{\mu+1-\ell}+\Span\{V_{\mu+2-\ell},...,V_{\ell+1},Y\},&\ell+1\le \mu\le 2\ell-1\\
\SS(\xi)^{2\ell}U_{\ell+1}&=&\star V_{\ell+1}+\Span\{Y\}.\end{array}$$
Thus $\SS(\xi)^{2\ell}\ne0$ for generic $\xi$. One shows similarly $\SS(\xi)^{2\ell+1}=0$
for every $\xi$ by:
$$\begin{array}{llll}
\SS(\xi)^\mu\BB&\subseteq&\Span\{T,U_2,...,U_{\ell+1-\mu},V_2,...,V_{\ell+1},Y\},&1\le
\mu\le\ell-1\\
\SS(\xi)^{\ell}\BB&\subseteq&\Span\{T,V_2,...,V_{\ell+1},Y\},\\
\SS(\xi)^{\mu}\BB&\subseteq&\Span\{V_{\mu+1-\ell},...,V_{\ell+1},Y\},&\ell+1\le \mu\le
2\ell\end{array}$$
and $\SS(\xi)^{2\ell+1}\BB =\{0\}$.
\end{proof}

We complete the proof of Theorem \ref{thm1.1} by considering the case $n=2\ell+2$ for
$\ell\ge1$. Let
$(x,u_1,u_2,...,u_{\ell+1},v_1,...,v_{\ell+1},y)$ be coordinates on
$\mathbb{R}^{2\ell+4}$ which define the
associated coordinate frame:

$$\BB:=\{X,U_1,...,U_{\ell+1},V_1,...,V_{\ell+1},Y\}=(\partial_x,\partial_{u_1},...,\partial_{u_{\ell+1}},
\partial_{v_1},...,\partial_{v_{\ell+1}},\partial_y\}.$$

\begin{lemma}\label{lem3.3} Let $\ell\ge1$. Define a metric $g_{2\ell+2}$ on
$\mathbb{R}^{2\ell+4}$ by setting:
\begin{eqnarray*}
&&g_{2\ell+2}(X,X)=-\textstyle\sum_{1\le a\le\ell}(u_a+v_a)u_{a+1}^2-\frac13u_1^3,\\
&&g_{2\ell+2}(X,Y)=1,\quad
g_{2\ell+2}(U_a,V_b)=\delta_{ab}.\end{eqnarray*}
Then $g_{2\ell+2}$ is a metric of signature $(\ell+2,\ell+2)$ and Szab\'o nilpotent of
order $2\ell+2$.\end{lemma}

\begin{proof} Let $2\le a\le\ell+1$. The non-zero components of $\nabla R_{g_{2\ell+2}}$
are:
$$\begin{array}{lllll}\nabla R_{g_{2\ell+2}}(X,U_1,U_1,X;U_1)&=&1,\\
\nabla R_{g_{2\ell+2}}(X,U_a,U_a,X;U_{a-1})
&=&\nabla R_{g_{2\ell+2}}(X,U_a,U_{a-1},X;U_a)&=&1,\\
\nabla R_{g_{2\ell+2}}(X,U_a,U_a,X;V_{a-1})
&=&\nabla R_{g_{2\ell+2}}(X,U_a,V_{a-1},X;U_a)&=&1.
\end{array}$$
We compute:
$$\begin{array}{lll}\SS(\xi)X&=&
    \star U_1+...+\star U_\ell+\star Y+\star V_1+...+\star V_{\ell+1},\\
\SS(\xi)U_1&=&\star Y+\star V_1+\star V_2,\\
\SS(\xi)U_a&=&
\star U_{a-1}+\star Y+\star V_{a-1}+\star V_a+\star
V_{a+1}\quad\text{for}
\quad2\le
a\le\ell,\\
\SS(\xi)U_{\ell+1}&=&\star U_\ell+\star Y+\star V_\ell
     +\star V_{\ell+1},\\
\SS(\xi)Y&=&0,\\
\SS(\xi)V_a&=&\star Y+\star V_{a+1}\quad\text{for}\quad 1\le
a\le\ell,\\
\SS(\xi)V_{\ell+1}&=&0.\end{array}$$
We may then show $\SS(\xi)^{2\ell+1}$ is generically non-zero by computing:
$$\begin{array}{llll}
S(\xi)^{\mu}U_{\ell+1}&=&\star U_{\ell+1-\mu}+\Span\{Y,V_1,...,V_{\ell+1}\},&1\le
\mu\le\ell\\
S(\xi)^{\mu}U_{\ell+1}&=&\star
V_{\mu-\ell}+\Span\{Y,V_{\mu+1-\ell},...,V_{\ell+1}\},&\ell+1\le \mu\le2\ell\\
S(\xi)^{2\ell+1}U_{\ell+1}&=&\star V_{\ell+1}+\Span\{Y\}.\end{array}$$
A similar argument shows $\SS(\xi)^{2\ell+2}=0$ for all $\xi$. We can write
$$\begin{array}{llll}
S(\xi)^{\mu}\BB&\subseteq&\Span\{U_1,...,U_{\ell+1-\mu},V_1,...,V_{\ell+1},Y\},&1\le
\mu\le\ell,\\
S(\xi)^{\mu}\BB&\subseteq&\Span\{V_{\mu-\ell},...,V_{\ell+1},Y\},&\ell+1\le
\mu\le2\ell+1\end{array}$$
and $S(\xi)^{2\ell+2}\BB = \{0\}$.
\end{proof}

\begin{remark}\label{rmk3.4} \rm One can also consider the purely pointwise question. We
shall say that $(M,g)$ is Szab\'o
nilpotent of order
$n$ at
$P\in M$ if
$\SS(x)^n=0$ for all $x\in T_PM$ and if $\SS(x_0)^{n-1}\ne0$ for some $x_0\in T_PM$.
Throughout Section \ref{Sect3}, we
considered cubic functions to ensure that $\nabla R$ was constant on the coordinate
frames; thus the point in question played
no role. However, had we replaced $u^3$ by $u^4$, $tu^2$ by $tu^3$, $u_au_{a+1}^2$ by
$u_au_{a+1}^3$, and
$v_au_{a+1}^2$ by $v_au_{a+1}^3$, then we would have constructed metrics $g_n$ which were
Szab\'o nilpotent of order $n$ on
$T_P\mathbb{R}^{n+2}$ for generic points $P\in\mathbb{R}^{n+2}$, but where $\nabla R$
vanishes at the origin
$0\in\mathbb{R}^{n+2}$. Since the order of nilpotency would vary with the point of the
manifold, these metrics clearly are not homogeneous.
\end{remark}

 \section{Nilpotent Osserman manifolds}\label{Sect4}
In Section \ref{Sect3}, we used cubic expressions to define our metrics to ensure the
tensors $R_{ijkl;n}$ were constant on
the coordinate frame. To discuss the Jacobi operator, we use the corresponding quadratic
polynomials. We adopt the
notation of Section \ref{Sect3} to define metrics:
\begin{eqnarray*}
&&\textstyle \tilde g_2(X,X)=-u^2,\ 
  \tilde g_2(X,Y)=1,\ \tilde g_2(U,V)=1,\\
&&\tilde g_3(X,X)=\textstyle-2tu-u^2,\ \tilde g_3(T,T)=1,\ 
  \tilde g_3(U,V)=1, \ \tilde g_3(X,Y)=1,\\
&&\tilde g_{2\ell+1}(X,X)=-2
   tu_2-u_2^2-\textstyle\sum_{2\le
   a\le\ell}\{2(u_a+v_a)u_{a+1}+u_{a+1}^2\},\\ 
&&\phantom{bork}\tilde g_{2\ell+1}(X,Y)=1,\
  \tilde g_{2\ell+1}(T,T)=1,\ 
  \tilde g_{2\ell+1}(U_a,V_b)=\delta_{uv},\quad(\ell\ge2)\\
&&\tilde g_{2\ell+2}(X,X)=-\textstyle\sum_{1\le
a\le\ell}\{2(u_a+v_a)u_{a+1}+u_{a+1}^2\}-u_1^2,\\ 
&&\phantom{bork} \tilde g_{2\ell+2}(X,Y)=1,\quad
\tilde g_{2\ell+2}(U_a,V_b)=\delta_{ab}\quad\qquad\qquad\qquad\qquad
(\ell\ge1).
\end{eqnarray*}
\begin{lemma}\label{lem4.1}\ \begin{enumerate}
\item $\tilde g_2$ has signature $(2,2)$ and is Osserman nilpotent of order $2$.
\item $\tilde g_3$ has signature $(2,3)$ and is Osserman nilpotent of
order
$3$.
\item $\tilde g_{2\ell+1}$ has signature $(\ell+1,\ell+2)$ and is Osserman
nilpotent of order
$2\ell+1$.
\item $\tilde g_{2\ell+2}$ has signature $(\ell+2,\ell+2)$ and is Osserman
nilpotent of order $2\ell+2$.
\end{enumerate}
\end{lemma}

\begin{proof} By Lemma \ref{lem2.1}, the non-zero components of $R_{\tilde g_2}$ are
\begin{equation}\label{eqn4.a}R_{\tilde g_2}(X,U,U,X)=1.\end{equation}
Assertion (1) now follows since:
$$J_{\tilde g_2}(\xi)X=\star Y+\star V,\quad
J_{\tilde g_2}(\xi)U=\star Y+\star V,\quad
J_{\tilde g_2}(\xi)Y=J(\xi)V=0$$
where $\star$ denotes suitably chosen quadratic polynomials in the components of $\xi$
which are non-zero for generic
$\xi$.

Similarly, the only non-zero component of $\nabla R_{\tilde g_3}$ are
\begin{equation}\label{eqn4.b}R_{\tilde g_3}(X,U,U,X)=1\quad\text{and}\quad
R_{\tilde g_3}(X,U,T,X)=1.\end{equation}
Assertion (2) now follows since:
$$\begin{array}{ll}
J_{\tilde g_3}(\xi)X=\star T+\star Y+\star V,&
J_{\tilde g_3}(\xi)Y=0,\\
J_{\tilde g_3}(\xi)T=\star Y+\star V,\\
J_{\tilde g_3}(\xi)U=\star T+\star Y+\star V,&
J_{\tilde g_3}(\xi)V=0.
\end{array}$$

We take $\ell\ge 2$ to prove assertion (3). Let $2\le a\le\ell$. The non-zero components
of
$R_{2\ell+1}$ are:
\begin{equation}\label{eqn4.c}\begin{array}{l}
1=R_{\tilde g_{2\ell+1}}(X,U_2,U_2,X)=
 R_{\tilde g_{2\ell+1}}(X,T,U_2,X)\\
\phantom{1}= R_{\tilde g_{2\ell+1}}(X,U_{a+1},U_{a+1},X)=
 R_{\tilde g_{2\ell+1}}(X,U_{a+1},U_a,X)\\
\phantom{1}=R_{\tilde g_{2\ell+1}}(X,U_{a+1},V_a,X).
\end{array}\end{equation}
Assertion (3) follows from the same argument as that used to prove Lemma \ref{lem3.2} as:
$$\begin{array}{l}
J(\xi)X\in\Span\{Y,T,U_2,...,U_\ell,
        V_2,...,V_{\ell+1}\},\\
J(\xi)Y=0,\\
J(\xi)T=\star Y+\star V_2,\\
J(\xi)U_2=\star T+\star Y+\star V_2+\star V_3,\\
J(\xi)U_a=\star U_{a-1}+\star Y+\star V_{a-1}+\star V_a+\star V_{a+1}\quad
\text{for}\quad3\le a\le\ell,\\
J(\xi)U_{\ell+1}=\star U_\ell+\star Y+\star V_\ell+\star V_{\ell+1},\\
J(\xi)V_a=\star Y+\star V_{a+1}\quad\text{for}\quad2\le a\le\ell,\\
J(\xi)V_{\ell+1}=0
\end{array}
$$

To prove assertion (4), we take $\ell\ge1$. Let $2\le a\le\ell+1$. The non-zero
components of
$R_{\tilde g_{2\ell+2}}$ are:
\begin{equation}\label{eqn4.d}\begin{array}{l}
 1=R_{\tilde g_{2\ell+2}}(X,U_1,U_1,X)=
 R_{\tilde g_{2\ell+2}}(X,U_a,U_a,X)\\\phantom{1}=
R_{\tilde g_{2\ell+2}}(X,U_a,U_{a-1},X)
=
R_{\tilde g_{2\ell+2}}(X,U_a,V_{a-1},X).
\end{array}\end{equation}
We may then compute:
$$\begin{array}{lll}J(\xi)X&=&
    \star U_1+...+\star U_\ell+\star Y+\star V_1+...+\star V_{\ell+1},\\
J(\xi)U_1&=&\star Y+\star V_1+\star V_2,\\
J(\xi)U_a&=&
\star U_{a-1}+\star Y+\star V_{a-1}+\star V_a+\star
V_{a+1}\quad\text{for}
\quad2\le
a\le\ell,\\
J(\xi)U_{\ell+1}&=&\star U_\ell+\star Y+\star V_\ell
     +\star V_{\ell+1},\\
J(\xi)Y&=&0,\\
J(\xi)V_a&=&\star Y+\star V_{a+1}\quad\text{for}\quad 1\le a\le\ell,\\
J(\xi)V_{\ell+1}&=&0.\end{array}$$
Assertion (4) now follows from the argument used to establish Lemma \ref{lem3.3}.
\end{proof}

\begin{remark}\label{rmk4.2} \rm Again, one can consider pointwise questions. We shall
say that $(M,g)$ is Ossersman nilpotent
of order
$n$ at
$P\in M$ if
$J(x)^n=0$ for all $x\in T_PM$ and if $J(x_0)^{n-1}\ne0$ for some $x_0\in T_PM$. By
replacing $u^2$ by $u^3$, $tu$ by
$tu^2$,
$u_au_{a+1}$ by $u_au_{a+1}^2$, and
$v_au_{a+1}$ by $v_au_{a+1}^2$, we could construct metrics $\tilde g_n$ on
$\mathbb{R}^{n+2}$ which are Osserman of order $n$
on $T_P\mathbb{R}^{n+2}$ for generic points $P\in\mathbb{R}^{n+2}$, but where $R$
vanishes at the origin
$0\in\mathbb{R}^{n+2}$. This gives rise to metrics where the order of nilpotency varies
with the point of the
manifold; such examples, clearly, are neither symmetric nor homogeneous.
\end{remark}

\begin{remark}\label{rmk4.3}\rm Stanilov and Videv 
\cite{StVi98} defined a higher order analogue of the Jacobi operator in the Riemannian
setting which was subsequently extended
to arbitrary signature. Let
$\Gr_{r,s}(M,g)$ be the Grassmannian bundle of all non-degenerate subspaces of $TM$ of
signature $(r,s)$. We assume $0\le r\le
p$,
$0\le s\le q$, and $0<r+s<p+q$ to ensure $\Gr_{r,s}(M,g)$ is non-empty and does not
consist of a single point; such a pair
$(r,s)$ will be said to be {\it admissible}. Let 
$\BB=\{e_1^+,...,e_r^+,e_1^-,...,e_s^-\}$ be an orthonormal basis for
$\pi\in\Gr_{r,s}(M,g)$. Then
$$J(\pi):=J(e_1^+)+...+J(e_r^+)-J(e_1^-)-...-J(e_s^-)$$
is independent $\BB$ and depends only on $\pi$. Following Stanilov, one says that $(M,g)$
is Osserman of type $(r,s)$ if the
eigenvalues of $J(\pi)$ are constant on $\Gr_{r,s}(M,g)$. Let $J_n$ be defined by the
metric $\tilde g_n$ defined in
Lemma
\ref{lem4.1}. The discussion given above then implies $J_n(\pi)^n=0$ for all $\pi$ and
thus $(\mathbb{R}^{n+2},\tilde g_n)$ is
Osserman of type $(r,s)$ for all admissible $(r,s)$. We refer to the discussion
in \cite{BCG01,GIZ02a} for other examples of higher order Osserman manifolds.
\end{remark}

\section{Ivanov-Petrova manifolds}\label{Sect5}
\begin{lemma}\label{lem5.1} The pseudo-Riemannian manifold $(\mathbb{R}^{n+2},\tilde
g_n)$ defined in {\rm Lemma
\ref{lem4.1}} is nilpotent Ivanov-Petrova of order $2$ if $n=2$ and nilpotent
Ivanova-Petrova of order $3$ if $n\ge3$.
\end{lemma}

\begin{proof} Suppose first $n=2$. We use equation (\ref{eqn4.a}) to see:
\begin{eqnarray*}
\mathcal{R}_{\tilde g_2}(\pi)X=\star V,\ \mathcal{R}_{\tilde g_2}(\pi)U=\star Y,\ 
\mathcal{R}_{\tilde g_2}(\pi)V=\RR_{\tilde g_2}(\pi)Y=0,
\end{eqnarray*}
where $\star$ are suitably chosen quadratic polynomials in the components of the generating vectors of $\pi=\Span\{f_1,f_2\}$ which
are non-zero for generic
$f_i$. Thus $\mathcal{R}_{\tilde g_2}(\pi)\ne0$ for generic $\pi$ while
$\mathcal{R}_{\tilde
g_2}(\pi)^2=0$ for all
$\pi$.

We use equations (\ref{eqn4.b}), (\ref{eqn4.c}), and (\ref{eqn4.d}) to compute
$\RR_{\tilde g_n}(\pi)Y=0$ and:
$$\begin{array}{ll}
\RR_{\tilde g_3}(\pi)X\in\Span\{ V, T\},&\RR_{\tilde g_3}(\pi)T\in\Span\{ Y\},\\
\RR_{\tilde g_3}(\pi)V=0,&\RR_{\tilde g_3}(\pi)U\in\Span\{ Y\},\\
\RR_{g_{2\ell+1}}(\pi)X\in\Span\{T,U_a,V_a\},&
\RR_{g_{2\ell+1}}(\pi)T\in\Span\{ Y\},\\
\RR_{g_{2\ell+1}}(\pi)U_a\in\Span\{Y\},&
\RR_{g_{2\ell+1}}(\pi)V_a\in\Span\{Y\},\\
\mathcal{R}_{\tilde g_{2\ell+2}}(\pi)X\in\Span\{U_a,V_a\},&
\mathcal{R}_{\tilde g_{2\ell+2}}(\pi)U_a\in\Span\{ Y\},\\
\mathcal{R}_{\tilde g_{2\ell+2}}(\pi)V_a\in\Span\{ Y\}.
\end{array}$$
This shows $\RR_{\tilde g_n}^3(\pi)=0$ $\forall$ $\pi$ and $\RR_{\tilde g_n}^2(\pi)\ne0$
for
generic $\pi$. \end{proof}

\section*{acknowledgments}
Research of P.G. partially supported by the NSF (USA) and MPI (Germany). 


\begin{thebibliography}{AAA}

\bibitem{BBG97} N. Bla\v zic, N. Bokan, and P. Gilkey, A
    note on Osserman Lorentzian manifolds, {\it Bull. London Math. Soc.},  {\bf 29},
       (1997), 227--230.

\bibitem{BBGR97} N. Bla\v zi\'c, N. Bokan, P. Gilkey and Z. Raki\'c,
     Pseudo-Riemannian Osserman manifolds, {\it Balkan J. Geom. Appl.}, {\bf 2}, (1997),
1--12.


\bibitem{BCG01} A. Bonome, P. Castro, E. Garcia-Rio, 
Generalized Osserman four-dimensional manifolds,
{\it Classical and Quantum Gravity}, {\bf 18} (2001), 4813--4822.

\bibitem{C88} Q.-S. Chi, A curvature characterization of certain locally
rank-one symmetric spaces, {\it J. Differential Geom.}, {\bf 28}, (1988), 187--202.

\bibitem{GKV97} E. Garc\'ia-Rio, D. Kupeli, and M. V\' azquez-Abal,
On a problem of Osserman in Lorentzian geometry, {\it Differential Geom. Appl.}, {\bf 7},
(1997),
85--100.

\bibitem{GVV98} E. Garc\'ia-Ri\'o, M. E. V\' azquez-Abal and
     R. V\' azquez-Lorenzo, Nonsymmetric Osserman pseudo-Riemannian manifolds,
     {\it Proc. Amer. Math. Soc.}, {\bf 126}, (1998), 2771--2778.

\bibitem{GVV02} E. Garci\'a-Ri\'o, D. Kupeli, and R. V\'azquez-Lorenzo, {\bf Osserman
Manifolds in Semi-Riemannian
Geometry}, Lecture notes in Mathematics 1777, Springer Verlag, Berlin, (2002), ISBN
3-540-43144-6.

\bibitem{G01} P. Gilkey, {\bf Geometric properties of natural operators defined by the
Riemann
curvature tensor}, World Scientific Publishing Co., (2001), ISBN 981-02-4752-4.

\bibitem{F01} B. Fiedler,  Determination of the structure of algebraic curvature tensors
by means of Young symmetrizers,
(preprint).

\bibitem{F02} B. Fiedler, On the symmetry classes of the first covariant derivatives of
tensor fields, (preprint).


\bibitem{GIS02} P. Gilkey, R. Ivanova, and I. Stavrov, 
Spacelike Jordan Szab\'o algebraic covariant curvature tensors in the higher
signature setting, preprint http://arXiv.org/abs/math.DG/0211089.

\bibitem{GIZ02} P. Gilkey, R. Ivanova, and T. Zhang,
Szab\'o Osserman IP Pseudo-Riemannian manifolds, \newline preprint 
http://arXiv.org/abs/math.DG/0205085

\bibitem{GIZ02a} ---, Higher order Jordan Osserman pseudo-Riemannian manifolds, \newline
preprint http://arXiv.org/abs/math.DG/0205269

\bibitem{GS02} P. Gilkey and I. Stavrov,  Curvature tensors
   whose Jacobi or Szab\'o operator is nilpotent
    on null vectors, {\it Bull. London Math. Soc.}, to appear;\newline preprint 
http://arXiv.org/abs/math.DG/0205074

\bibitem{KTV92} O. Kowalski, F. Tricerri and L. Vanhecke,
Curvature homogeneous Riemannian manifolds,
{\it J. Math. Pures Appl.}, {\bf 71} (1992), 471--501.

\bibitem{N02} Y. Nikolayevsky, Osserman Conjecture in dimension $n \ne 8, 16$,
preprint:\newline
http://arXiv.org/abs/math.DG/0204258.

\bibitem{Oss90} R. Osserman, Curvature in the eighties, {\it Amer. Math.
    Monthly}, {\bf97}, (1990) 731--756.

\bibitem{StVi98} G. Stanilov and V. Videv, Four dimensional
    pointwise Osserman manifolds, {\it Abh. Math. Sem. Univ. Hamburg}, {\bf 68}, (1998),
    1--6.

\bibitem{St03} I. Stavrov, Ph. D. Thesis, University of Oregon (2003).

\bibitem{Sz91}
Z. I. Szab\'o, A short topological proof for the symmetry of $2$ point homogeneous
spaces, {\it Invent. Math.}, {\bf 106}, (1991), 61--64.

\end{thebibliography}
\end{document}